\documentclass[12pt]{article}
\usepackage{amsmath,amssymb,amsthm,epsfig,psfrag}

\newtheorem{theorem}{Theorem}

\newtheorem{proposition}[theorem]{Proposition}

\newtheorem{lemma}[theorem]{Lemma}

\theoremstyle{definition}

\def\cal{\mathcal}

\def\Bbb{\mathbb}

\def\Z{\Bbb{Z}}

\def\ni{\noindent}

\def\F+L{\hbox{$\textup{F}\!_+\textup{L}$}}

\def\ssm{\smallsetminus}

\def\ms{\medskip}

\def\onto{{\kern3pt\to\kern-8pt\to\kern3pt}}

\def\<{\langle}
\def\>{\rangle}
\def\|{{\ |\ }}

\def\II{\cal I}

 \def\AA{\cal A}

 \def\RR{\cal R}

\def\LL{\cal L}
 \def\PP{\cal P}

\newcommand{\set}[1]{\left\{#1\right\}}

\renewcommand{\ni}{\noindent}

\renewcommand{\ms}{\medskip}
\newcommand{\bs}{\bigskip}
\def\*{^{\star}}
\newcommand{\Proof}{\ni \emph{Proof. }}

\def\qed{{\ifhmode\unskip\nobreak\hfil\penalty50 \hskip1em \else\nobreak\fi
   \mbox{}\nobreak\hfil \rule{1ex}{1ex}
   \parfillskip=0pt \finalhyphendemerits=0 \par}}

\begin{document}

\title{A finitely presented group \\ with unbounded dead-end depth}

\author{Sean Cleary\thanks{Support from PSC-CUNY grant \#65752  is gratefully acknowledged.}
 and Tim R.\ Riley\thanks{Support from NSF grant  0404767 is gratefully acknowledged.}}

\date{June 22, 2004 \\ 
Revised August 25, 2004 and April 5, 2007}
\maketitle

\begin{abstract}
\ni The \emph{dead-end depth} of an element $g$ of a group $G$, with respect to a generating set 
$\AA$, is the distance from $g$ to the complement of the radius $d_{\AA}(1,g)$ closed ball, in the 
word metric $d_{\AA}$ defined with respect to $\AA$.   We exhibit a finitely presented group $G$ with a finite generating set with respect to which there is no upper bound on the dead-end depth of elements. \\
\footnotesize{\ni \textbf{2000 Mathematics Subject
Classification: 20F65}  \\ \ni \emph{Key words and phrases: dead-end depth, lamplighter}}

\bs

\ni\emph{The authors regret that the published version of this article \textup{(}Proc.\ Amer.\ Math.\ Soc., 134\textup{(}2\textup{)}, pages 343-349, 2006\textup{)} contains a significant error concerning the model for $G$ described in Section~2.  We are grateful to J\"org~Lehnert for pointing out our mistake.  In this corrected version, that model has been overhauled, and that has necessitated a number of changes in the subsequent arguments.}
\end{abstract}

\section{Introduction} \label{intro}

We explore the behavior of geodesic rays and the geometry of balls in Cayley graphs of
finitely generated groups.
An element $g$ in a group $G$ is defined by Bogopol$'$ski\u{\i}~\cite{bogop} to be a  {\em dead end} with respect to a generating set $\AA$ (which will always be finite in this article) if it is not adjacent to an element further from the identity; that is, if a geodesic ray in the Cayley graph of $(G, \AA)$ from the identity to $g$ cannot be extended beyond $g$.  Dead-end elements occur in a variety of settings.  Bogopol$'$ski\u{\i} showed that $\textup{SL}_2(\Z)$, presented by $\langle x,y \mid x^4, y^6, x^{-2}y^3 \rangle$, has dead ends with respect to $\set{x,y}$ and that $\langle x,y \mid x^3, y^3, (xy)^k \rangle$, with $k \geq 3$ has dead ends with respect to $\set{x,y}$  but not $\set{x,y, xy}$.  Fordham
\cite{blake:diss} found dead ends in  Thompson's group $F$, presented by $\langle x_0, x_1 \mid [x_0{x_1}^{-1}, {x_0}^{-1}x_1x_0], [x_0{x_1}^{-1}, {x_0}^{-2}x_1{x_0}^2]\rangle$, with respect
to the standard finite generating set $\set{x_0,x_1}$.  In Lemma 4.19 of \cite{Champ} Champetier  shows that presentations $\langle \AA \mid \RR \rangle$ satisfying the $C'(1/6)$ small cancellation condition have no dead ends with respect to $\AA$.  In general, generating sets can be contrived with respect to which there are dead ends; for example, $\Z$ has dead ends with respect to the generating set $\{a^2, a^3\}$.  Many of these examples are discussed in IV.A.13,14 of de~la~Harpe \cite{delaharpe}, where there is also an exercise due to Valette concerning further examples.

Dead ends differ in their severity in the following sense.  We define the {\em depth}, with respect
to a finite generating set $\AA$, of an element $g$ in an infinite group $G$ to be the distance in the word metric $d_{\AA}$ between $g$ and the complement  in $G$ of the closed ball $B_g$ of radius $d_{\AA}(1,g)$ centered at $1$.
So $g \in G$ is a dead end when its depth is at least $2$. Cleary and Taback \cite{ctcomb} showed that, with respect to $\set{x_0,x_1}$, the dead-end elements in Thompson's group $F$ are all of depth $3$.  

Cleary and Taback \cite{deadlamp} exhibited wreath products, such as the lamplighter group $\Z_2 \wr \Z$, with unbounded dead-end length with respect to certain standard generating sets.  Independently, Erschler observed that $\Z_2 \wr \Z$ provides an example resolving a closely related question of Bowditch  (Question 8.4 in Bestvina's problem list \cite{mladenlist}): let $\Gamma$ be the Cayley graph of an infinite finitely generated group; does there exist $K > $0 such that for all $R > 0$ and all vertices $v \in \Gamma - \textit{B}(R)$ there is an infinite ray from $v$ to $\infty$ which does not enter $\textit{B}(R - K)$?  

However, these wreath product examples are not finitely presentable.  
So a natural question, asked by Bogopol$'$ski\u{\i}~\cite{bogop}, is whether, given a finitely presentable group and a finite generating set $\AA$, there is an upper bound on dead-end depth with respect to $\AA$.  Bogopol$'$ski\u{\i}~\cite{bogop} showed that there is always such a bound in the case of hyperbolic groups.  We answer this question in the negative with the main result of this paper:

\begin{theorem}  \label{mainthm}
There is a finitely presentable group $G$ that has a finite
generating set $\AA$ with respect to which $G$ has unbounded dead-end depth.  We  take $G$ to be the group presented by $$\PP:=
\langle \ a,s,t \ \mid \ a^2=1, \ [a,a^t]=1, \ [s,t]=1, \ a^s=aa^t \
\rangle$$ and $\AA$ to be the  generating set $$\set{a,s,t,at,ta,ata,as,sa,asa}.$$
\end{theorem}

We denote the commutator $a^{-1} b^{-1} a b$ by $[a,b]$ and 
the conjugate $b^{-1}ab$ by $a^b$.

The lamplighter group $\Z_2 \wr \Z$, which has presentation
$$\langle \ a,t \ \mid \ a^2, [a,a^{t^i}] , \ \forall i \in \Z \
\rangle,$$ is a subgroup of $G$.  Removing the defining relation
$a^2$ from $\PP$ gives Baumslag's remarkable example
\cite{gilbertmeta}  of a finitely presented metabelian group
containing $\Z \wr \Z$ and thus  a free abelian subgroup of
infinite rank.

The reference to a specific finite generating set in the theorem
is important as the following issues are unresolved.\footnote{\emph{Remark added August 2004. } It
has  since been shown \cite{RW} that $ \langle \ a, t, u  \mid a^2=1, \ [t,u]=1, \ a^t=a^u\,; \ \forall i \in \Z, \ [a,a^{t^i}] \ \rangle$ has unbounded depth dead ends with respect to one finite generating set, but 
only a single dead end with respect to another.}  
Is the property of having unbounded dead-end depth an invariant of
finitely generated groups?  That is, does it depend on the choice
of finite generating set?  Moreover, is this property a
quasi-isometry invariant?    Indeed, does the group $G$ defined
above have unbounded dead-end depth with respect to all finite
generating sets?

\ms \ni \emph{Acknowledgment.}  We began discussing the ideas in this paper at the 2004 Cornell Topology Festival.  We are grateful to the organizers for their hospitality.  We would also like to thank Oleg~Bogopol$'$ski\u{\i} for an account of the history of dead-end depth and the anonymous referee for a careful reading.

\section{The lamplighter grid model for $G$} \label{grid model}

The group $\Z_2 \wr \Z$ was named the \emph{lamplighter} group by Cannon (see Parry~\cite{parrynp}) 
on account of the following faithful,
transitive action on $\PP_{\!\!\textit{fin}}(\Z) \times
\Z$. Here, $\PP_{\!\!\textit{fin}}(\Z)$ denotes the set of finite
subsets of $\Z$; that is, finite configurations illuminated among a string of lamps indexed by
$\Z$. The second factor denotes the location of a
\emph{lamplighter} (or a cursor) in the string of lamps.  The action of the generator 
$t$ of $\langle  a,t   \mid   a^2, [a,a^{t^i}] , \ \forall i \in \Z  
\rangle$ is to increment the location of the lamplighter by one, and the action of $a$ is to toggle the lamp at the current location of the lamplighter between on
and off.

We will develop a lamplighter model for $G$ and describe a faithful, transitive left
action of $G$ on $\PP_{\!\!\textit{fin}}(\mathbb{L}) \times \Z^2$, where $\PP_{\!\!\textit{fin}}(\mathbb{L})$ denoted the set of finite subsets of a countable set $\mathbb{L}$ that we define below.
This will be be useful for gaining insights into the metric properties of $G$.

\begin{figure}[ht]
\psfrag{a}{$a$}%
\psfrag{s}{$s$}%
\psfrag{t}{$t$}%
\centerline{\epsfig{file=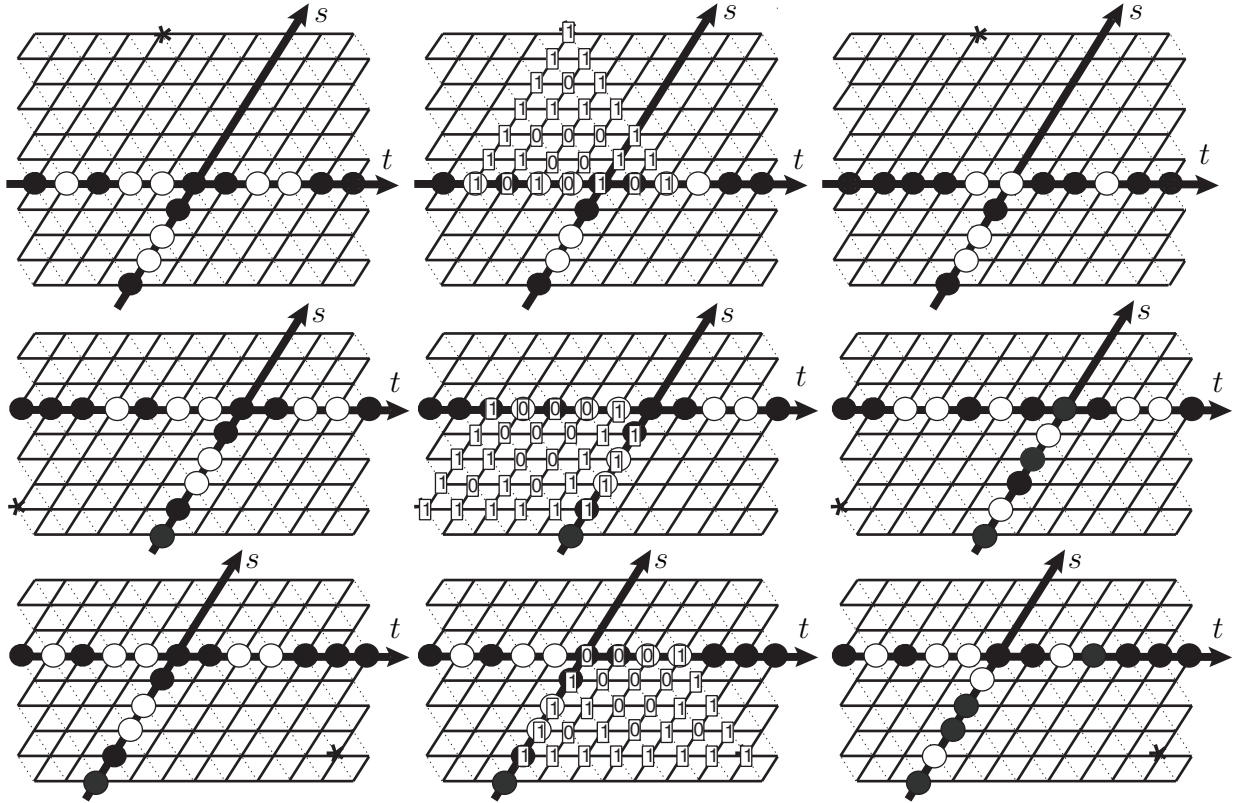}}
 \caption{Examples of $a$ acting.  The first column of diagrams shows $hg ( \emptyset, (0,0) )$, where $g=s^3as^{-1}as^{-2}t^4at^{-2}at^{-1}at^{-3}at^{-1}at^{-1}at^3$ and  $h$ is $s^{6}t^{-4}$,  $s^{-4}t^{-4}$,  and $s^{-4}t^{7}$, respectively.  The third column show the corresponding $a hg( \emptyset, (0,0) )$.  (Off and on lamps are represented  by circles with black and white interiors, respectively.)} \label{lamp-lighting grid}
\end{figure}

In our model, a lamplighter moves  among the
lattice points of the infinite rhombic grid illustrated in
Figure~\ref{lamp-lighting grid}.  We will refer to  the union of the $t$-axis and the portion of the $s$-axis that is below the $t$-axis as the \emph{lampstand}.  Let  $\mathbb{L}$ be the set of lattice points on the lampstand; these are the locations of the lamps in our model.   An element of
$\PP_{\!\!\textit{fin}}(\mathbb{L})$ denotes a finite configuration of illuminated lamps.  

The actions of $s$ and $t$ are to move the
lamplighter one unit in the $s$- and $t$-directions, respectively, in the rhombic grid.  The rhombic grid is
subdivided into a triangular grid by inserting a negatively sloped diagonal into each rhombus (the dashed lines in Figure~\ref{lamp-lighting grid}).  At every lattice point there is a \emph{button}; the action of $a$ is to \emph{press} the button at the location of the lamplighter, and this has the following effect.   If the button is on the lampstand, then it toggles the lamp at its location.
If the button is off the lampstand, a \emph{signal is set off} that propagates and bifurcates  in the triangular  grid towards the lampstand and toggles  finitely many lamps as follows.   

When the button is at a lattice point above the $t$-axis, the 
signal propagates downwards in the triangular grid along the sloped grid lines.  
At each vertex en route it splits into two signals, one advancing along the positively sloped diagonal below and one along the negatively sloped diagonal below.  The signals stop when they hit the $t$-axis, and each lamp on the $t$-axis switches between on
and off once for every signal it receives.   The manner in which these signals split as they propagate towards the lampstand
leads to a connection with Pascal's triangle modulo $2$.
When the lamplighter is above the $t$-axis, we can understand the action of $a$ as follows.  Suppose the current location of the lamplighter has $s,t$--coordinates $(p,q)$.  Then for $r \in \set{0,1, \ldots, p}$, the lamp at position $t^{q+ r}$ is toggled
between on and off when there is a $1$ in the $r$-th entry of row
$p$ of Pascal's triangle mod $2$.   See the top row of diagrams in Figure~\ref{lamp-lighting grid} for an example.  

When the button is at a lattice point below the $t$-axis and to the left of the $s$-axis, the signals propagate similarly through the triangular grid towards the lampstand, but move in the horizontal direction and in the $s$-direction.  The signals stop when they hit the lampstand.       
  There is again a connection with Pascal's triangle modulo $2$ but this time it is rotated by $2\pi/3$.  The pattern is interrupted at the $s$- and $t$-axes because signals stop when they hit the lampstand.  The middle  row of diagrams in Figure~\ref{lamp-lighting grid} shows an example.  

When the button is at a lattice point below the $t$-axis and to the right of the $s$-axis, the signals propagate  in the horizontal direction and in the negatively sloped diagonal direction towards the lampstand, where they again stop.  Pascal's triangle modulo $2$ rotated by $4\pi/3$ from its standard orientation (and interrupted on the $t$-axis) can be used in describing the action, as illustrated in the bottom row of diagrams of 
Figure~\ref{lamp-lighting grid}.

To verify that we have a well-defined action of $G$  on $\PP_{\!\!\textit{fin}}(\mathbb{L}) \times \Z^2$, it suffices to check that $a^2$, $[a,a^t]$, and $[s,t]$ all act
trivially and that the actions of  $aa^t$ and $a^s$ agree.  This and the proof that the action is
transitive are left  to the reader.  

To show that the action is faithful we suppose $g \in G$ satisfies $g(\emptyset, (0,0)\,) = (\emptyset, (0,0)\,)$ and we will check that $g=1$ in $G$.  Let $w$ be a word representing $g$.   Consider the rhombic grid with the lights all off and the lamplighter at $(0,0)$.  Reading $w$ from right to left determines a path for the lamplighter starting and finishing at $(0,0)$, in the course of which  buttons are pressed: each letter $s^{\pm 1}, t^{\pm 1}$ in $w$ moves the lamplighter, and each $a^{\pm 1}$ \emph{presses a button} at the lamplighter's location, toggling some of the lamps.  
Use the relations $a^2=1$ and $aa^t=a^s$ to change $w$ to a word $w'$ such that $w=w'$ in $G$, and there are no letters $a^{-1}$ in $w'$ and no subwords $a^2$, and along the path determined by $w'$, the only buttons that are pressed are on the lampstand.  Then, by  applying the relation $[s,t]=1$ to $w'$ and removing inverse pairs, obtain a freely reduced word $w''$ such that $w' = w''$ in $G$ and the path determined by $w''$ does not leave the lampstand.  Then by inserting inverse pairs of letters, $w''$ can be changed to a word $w'''$ whose path returns to the origin after each button-press --- that is, $w'''$ is a concatenation of words of the form $a^{s^i}$  (with $i <0$) and  $a^{t^j}$ (with $j \in \Z$).  Moreover, as $w'''$ leaves no lamps illuminated,   each $a^{s^i}$ and $a^{t^j}$ occurs an even number of times.  But, as we show below, all $a^{s^i}$ and $a^{t^j}$ represent commuting elements of $G$.  We deduce that $w'''$, and therefore $w$, represents $1$ in the group.

To see that $[a^{t^i}, a^{t^j}] =1$ in $G$ for all $i,j$, it suffices to show that $[a, a^{t^i}]=1$ for all $i \geq 0$.  We will induct on $i$.  Assume $[a, a^{t^i}]=1$ for all $0 \leq i \leq k$.  Then $[a,a^{t^k}]^s=1$ and so  $[a a^t, a^{t^k} a^{t^{k+1}}]=1$ as $a^s=aa^t$.  With the exception of $a, a^{t^{k+1}}$, the induction hypothesis tells us that any two of $a, a^t, a^{t^k}, a^{t^{k+1}}$ commute.  So $[a , a^{t^{k+1}}] = [a a^t, a^{t^k} a^{t^{k+1}}]= 1$. 

To see that $[a^{s^{i}}, a^{t^j}] =1$ for all $i,0$ and all $j\in \Z$, notice that
$[a^{s^{i}}, a^{t^j}]^{s^{-i}} = [a, a^{t^j s^{-i}}] = [a, (a^{s^{-i}})^{t^j}]$ which is $1$ in $G$ because $(a^{s^{-i}})^{t^j}$ is in the subgroup $S$ of $G$ generated by $\set{a^{t^k}}_{k \in \Z}$. To see that $[a^{s^{i}}, a^{s^{j}}] =1$ for all $i,j<0$, it suffices to show that $[a, a^{s^l}] =1$ for all $l>0$ and that can be done similarly.

\section{Proof of the theorem}

Define maps $\II:G \to \PP_{\!\!\textit{fin}}(\mathbb{L})$ and $\LL: G
\to \Z^2$ by $$(\, \II(g), \LL(g)\,)=g(\emptyset, (0,0)\,).$$ 
So $\II$ gives the lamps illuminated and $\LL$ gives the
location of the lamplighter after the action of $g$ on the configuration in which no lamps are lit and the lamplighter is at the origin.   (In fact, $\LL$ is, in effect, the retraction $G \onto \mathbb{Z}^2 = \langle s, t \rangle$ that kills $a$.)

\begin{figure}[ht]
\psfrag{1}{$L_1$}%
\psfrag{2}{$L_2$}%
\psfrag{3}{$L_3$}%
\psfrag{a}{$a$}%
\psfrag{s}{$s$}%
\psfrag{t}{$t$}%
\centerline{\epsfig{file=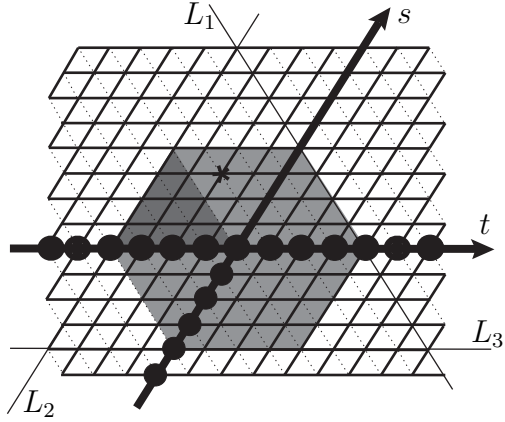}} \caption{Regions and lines in the lamplighting grid.} \label{grid2}
\end{figure}

Define $H_n$
to be the subset of $\Z^2$ of lattice points in (and on the
boundary of) the hexagonal region of the grid with corners at
$(\pm n,0), (0, \pm n), (n,-n), (-n,n)$.   (The  shaded region in Figure~\ref{grid2} is $H_4$.)

In the following proposition we determine the distance of various
group elements from the identity in the word metric with respect
to $$\AA = \set{a,s,t,at,ta,ata,as,sa,asa}.$$  The crucial
feature of this generating set is that the button at a vertex the
lamplighter is leaving or arriving at can be pressed with no
additional cost to word length.  Thus for $g \in G \ssm \set{a}$, the distance $d_{\AA}(g,1)$ is the 
length of the shortest path in the rhombic grid that starts at $(0,0)$, finishes at $\LL(g)$, and such that pressing some of the buttons at the vertices visited produces the configuration $\II(g)$ of illuminated bulbs.

\begin{proposition} \label{all the hard work}
    If  $g \in G$ is such that $\II(g) \subseteq H_n$ and $\LL(g) \in H_n$, then $d_{\AA}(1,g) \leq 6n$.  
\end{proposition}

\Proof
Define $(p,q) := \LL(g)$.  So $p$ and $q$ are the $s$-- and $t$--coordinates, respectively, of the position of the lamplighter.   We will describe a path from $(0,0)$ to $(p,q)$ and specify buttons to push en route that will illuminate $\II(g)$.

We will first address the case where $p \geq 0$.  Have the lamplighter begin by travelling a distance $2n$ first along the $s$-axis to $(-n,0)$ and then back to $(0,0)$, illuminating all bulbs in $\II(g)$ that are not on the $t$-axis as it goes.  

How to proceed next depends on whether  $(p,q)$ is in the set  $T_n$  of lattice points in, and on the
boundary of, the triangular region with corners $(0,0), (0,-n),
(n,-n)$. (In Figure~\ref{grid2} the region 
shaded dark gray contains $T_4$.)
Assume $(p,q) \in T_n$.    Have the lamplighter travel along the
$t$-axis to $(0,n)$, then along the $t$-axis to $(0,-n)$, then
along the $t$-axis to $(0,q)$, and then parallel to the $s$-axis
to $(p,q)$.  This involves traversing at most $4n$ edges.
 Following this path, the lamplighter visits all the lamps in $\II(g)$ that are on the $t$-axis
and so can illuminate them.  So $d_{\AA}(1,g) \leq 6n$.

Next, we assume $\LL(g) \in H_n \ssm T_n$.  Then
$\LL(g)$ is in, or on the boundary of, the rhombic region above
the axis (shaded medium gray in Figure~\ref{grid2}).  Have the lamplighter follow a path $\phi$, first along the $t$-axis to $(0,-n)$, then back along the $t$-axis to $(0,-p)$, then parallel
to the $s$-axis to $(p,-p)$, then parallel to the $t$-axis to
$(p,n-p)$, and then parallel to the $t$-axis back to to $(p,q)$. 
These five arcs have length $n$, $n-p$, $p$, $n$ and $n-p-q$,
respectively, totalling $4n-p-q$, which is less than $4n$ because
$p>-q$.

\begin{figure}[ht]
\psfrag{1}{$1$}%
\psfrag{0}{$0$}%
\centerline{\epsfig{file=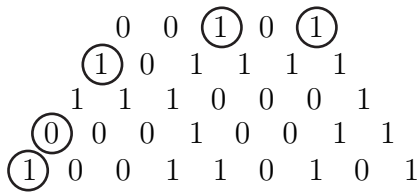}} \caption{A trapezoid of ones and zeros in which every entry other than those circled (called \emph{summits}) is the sum of the two entries above it mod $2$.} \label{trapezoid fig}
\end{figure}

We now show that pressing some combination of buttons on $\phi$
achieves the configuration of illuminated lamps on the $t$-axis required for $\II(g)$.  The
lamps in positions $-n$, $-n+1$, $\ldots$, $-p-1$ can be
illuminated as required when traversing the second arc of $\phi$.
 The lamps in positions $-p$, $-p+1$, $\ldots$, $n$ determine a
sequence of $n+p+1$ zeros and ones which make up the base of a
trapezoid of zeros and ones as analyzed in Lemma~\ref{trapezoid}, below. 
Overlaying this trapezoid on the grid, we find that its left and
top sides follow the third and fourth arcs of $\phi$.  
Pressing the buttons at the locations of the \emph{summits} (defined in
Lemma~\ref{trapezoid}) in the trapezoid illuminates the required lamps.  
This is because the trapezoid is constructed in such a way that each entry equals the total number of signals (see Section~\ref{grid model}) mod $2$ reaching its location when the buttons at the summits are all pressed.  In other words, the trapezoid is the overlay (mod 2) of a number of copies of Pascal's triangle, one for each summit; so, for example, the trapezoid in Figure~\ref{trapezoid fig} is the overlay of five triangles, with heights $1,2,4,5$ and $5$.  

We now turn to the case where $p<0$ and $q \leq 0$.  The lamplighter begins by travelling along the $t$-axis to $(0,n)$ and then back to $(0,0)$, illuminating all  the bulbs in $\II(g)$ that it traverses en route.  Next it moves along the $t$-axis to $(0,-n)$ and then parallel to the $s$-axis to $(-n,-n)$.  As in the argument above involving Lemma~\ref{trapezoid}, between $(0,-n)$ and $(-n,-n)$ the lamplighter can press a combination of buttons which put the lights on the $s$-axis into the configuration of $\II(g)$.  This will have some effect on the lights between $(-1,0)$ and $(-n,0)$.  But, when travelling from $(0,0)$ to $(0,-n)$, the lamplighter can toggle lights  both to counteract this effect and to implement the required configuration for those lights.  By this stage the lamplighter has travelled a distance of $4n$, is located at $(-n,-n)$, and $\II(g)$ has been achieved.  All that remains is for the lamplighter to move to $(p,q)$, which is within a further distance of $2n$ as $p<0$ and $q \leq 0$.  

Finally we address the case where $p<0$ and $q > 0$.  First the lamplighter travels along the $t$-axis to $(0,-n)$ and then back to $(0,0)$, illuminating all the bulbs as per $\II(g)$ en route.  Then the lamplighter moves along the $t$-axis to  and then parallel to the $t$-axis to $(-n, n)$.  Between $(0,n)$ and $(-n, n)$ buttons are pressed to illuminate the lamps on the $s$-axis as per $\II(g)$, and between $(0,0)$ and $(0,n)$ buttons are pressed to counteract the effect on the bulbs there and to achieve the configuration $\II(g)$.  The lamplighter then moves  to $(p,q)$ which is within a distance $2n$ from $(-n,n)$. \qed


\begin{lemma} \label{trapezoid}
Suppose $S$ is a sequence of $m$ zeros and ones.  Suppose $r \in \set{1, 2, \ldots, m}$. There is a trapezoid consisting of zeros and ones \emph{(}\emph{entries}\emph{)} satisfying the following.  The entries are arranged in $r$ rows with the bottom row $S$ and all the other rows containing one fewer entry than that below it \emph{(}as in Figure~\ref{trapezoid fig}\emph{)}.  Every entry is the sum of the two entries immediately above it mod 2, with possible exception of some entries in the top row and at the left ends of the rows.  We call these exceptional entries \emph{summits}.
\end{lemma}

\Proof We use induction on $r$.  When $r=1$ there is nothing to prove.  For the induction step, we assume the lowest $r$ rows have been constructed as per the hypotheses.  We add the next row one entry at a time starting from the right-hand side.  We select each entry so as to ensure that of the two entries immediately below it, that to the right has the property that it is the sum of the two entries immediately above it mod 2.  This ensures that only entries in the left side the trapezoid and in the $(r+1)$-st row can fail to equal the sum of the two entries immediately above them, mod 2.      
\qed

\begin{proposition} \label{lower bound}
If $g \in G$ satisfies $\LL(g)= (0,0)$ and $\II(g) =\set{(0,n), (0,-n), (-n,0)}$,  then 
$d_{\AA}(1,g) \geq 6n$.
\end{proposition}

\Proof
In order to light the lamps at $(0,n)$, $(0,-n)$ and $(-n,0)$, the lamplighter must visit (or cross) the straight lines $L_1$ through $(0,n)$ and $(n,0)$, $L_2$ through $(0,-n)$ and $(n,-n)$, and $L_3$ through $(-n,0)$ and $(-n,n)$, and then return to the origin.  (See Figure~\ref{grid2}.)  We suppose the lamplighter  follows a minimal length path that visits the three lines and returns to the origin.  We will show that this path has length at least $6n$.   

 Drawing $L_1, L_2, L_3$ in standard $\mathbb{R}^2$ with the $\ell_1$-metric, a symmetry becomes apparent on account of which we may as well assume that the lamplighter first visits $L_2$.  But then because we are using the $\ell_1$-metric, the path can be altered without increasing its length so that it first travels from the origin to $L_2$ along the $t$-axis.  If  the path visits $L_3$ before $L_1$, then it can be altered without increasing its length so that it does so by travelling a distance $n$ in the (negative)  $s$-direction --- one then easily checks that visiting $L_1$ and then returning to the origin costs at least $4n$ in additional length.  Suppose, on the other hand, the lamplighter proceeds next to $L_1$ after $L_2$.  Then the final portion of the lamplighters journey will be between $L_3$ and the origin and the path can be altered, with no change in length, so that it follows the $s$-axis from $(-n,0)$ to $(0,0)$.  One checks that a shortest path from $(0,-n)$ to  $(-n,0)$ via $L_1$ has length at least $4n$.   

So in each case the length of the path is at least $6n$ and so $d_{\AA}(1,g) \geq 6n$,  as claimed. 
\qed

\ni \emph{Proof of Theorem \ref{mainthm}.}
Define $$g_n \ :=  \ s^nas^{-n}t^n a t^{-2n} a t^n \ = \  s^{n-1} (sa) s^{-n}t^{n-1} (ta) t^{-2n} (at) t^{n-1}.$$  Then $d_{\AA}(1,g_n) \leq 6n$ and so $d_{\AA}(1,g_n) = 6n$ by Proposition~\ref{lower bound} because $g_n( \emptyset, (0,0) ) = ( \set{(-n,0),(0,n),(0,-n) }, (0,0) )$.  By Proposition~\ref{all the hard work}, if $w$ is a word representing an element $g$ with $d_{\AA}(1,g) > 6n$, then the lamplighter's path determined by $w$ (reading right to left and with the lamplighter initially at the origin) must leave $H_n$.  But
$\LL(g_n)=(0,0)$ and $H_n$ contains the closed ball of radius $n$ about $(0,0)$ in the rhombic grid.  So if $v_n$ is a word representing $g_n$ and $u_nv_n$ a word representing an element $h_n$ with $d_{\AA}(1, h_n) > 6n$, then the length of $u_n$ is more than $n$.  So, with respect to $d_{\AA}$, the distance from $g_n$ to the complement of the radius $d_{\AA}(1,g_n)$ ball is more than $n$.   Thus, with respect to $\AA$, the dead-end depth of $g_n$ is at least $n$ and $G$ has unbounded dead-end depth.
\qed

\bibliographystyle{plain}
\def\cprime{$'$}

\small{ \ni \textsc{Sean Cleary} \rule{0mm}{6mm} \\ Department of Mathematics, The City College of New York, 
City University of New York, New York, NY 10031,  USA  \\ \texttt{cleary@sci.ccny.cuny.edu, \
http:/\!/www.sci.ccny.cuny.edu/$\sim$cleary/ }

\ni  \textsc{Tim R.\ Riley} \rule{0mm}{6mm} \\  Department of Mathematics, 310 Malott Hall, Cornell University, Ithaca, NY 14853-4201, USA \\ \texttt{tim.riley@math.cornell.edu, \
http:/\!/www.math.cornell.edu/$\sim$riley/ } }

\end{document}